\documentclass{tbimc}
\usepackage{latexsym}
\usepackage{amsmath}
\usepackage{amssymb}

\newcommand{\eps}{\varepsilon}

\newcommand{\bi}{\begin{itemize}}
\newcommand{\ei}{\end{itemize}}
\newcommand{\beq}{\begin{equation}}
\newcommand{\eeq}{\end{equation}}
\newcommand{\abso}[1]{\left| #1 \right|}

\newcommand{\fin}{\begin{flushright} $\Box $ \end{flushright}}
\newcommand{\vsp}{\vspace{10mm}}

\newcommand{\matR}{\mathbb{R}}

\newtheorem{thm}{THEOREM}

\newtheorem{lm}{LEMMA}
\newtheorem{cor}{COROLLARY}

\begin{document}

\selectlanguage{english}

\title{L\`{e}vy approximation of impulsive recurrent process with Markov switching.}

\author{V. S. Koroliuk}
\address{}
\curraddr{Institute of Mathematics,\\ Ukrainian National Academy of
Science, Kiev, Ukraine} \email{korol@imath.kiev.ua}
\thanks{The authors thank University of Bielefeld
for hospitality and financial support by DFG project 436 UKR
113/80/04-07.}

\author{N. Limnios}
\address{}
\curraddr{Laboratoire de Math\'ematiques Appliqu\'ees,\\
Universit\'e de Technologie de Compi\`egne, France}
\email{nlimnios@dma.utc.fr}
\thanks{}

\author{I.V. Samoilenko}
\address{}
\curraddr{Institute of Mathematics,\\ Ukrainian National Academy of
Science, Kiev, Ukraine} \email{isamoil@imath.kiev.ua}
\thanks{}

\subjclass[2000]{Primary 60J55, 60B10, 60F17, 60K10; Secondary
60G46, 60G60.}
\date{09/02/2009}
\dedicatory{} \keywords{L\'{e}vy approximation, semimartingale,
Markov process, impulsive recurrent process, piecewise deterministic
Markov process, weak convergence, singular perturbation.}

\begin{abstract}
In this paper, the weak convergence of impulsive recurrent process
with Markov switching in the scheme of L\'{e}vy approximation is
proved. For the relative compactness, a method proposed by R.
Liptser for semimartingales is used with a modification, where we
apply a solution of a singular perturbation problem instead of an
ergodic theorem.
\end{abstract}

\maketitle

\section{Introduction}
L\'{e}vy approximation is still an active area of research in
several theoretical and applied directions. Since L\'{e}vy processes
are now standard, L\'{e}vy approximation is quite useful for
analyzing complex systems (see, e.g. \cite{ber, sat}). Moreover they
are involved in many applications, e.g., risk theory, finance,
queueing, physics, etc. For a background on L\'{e}vy process see,
e.g. \cite{ber, sat, gisk}.

In particular in \cite{korlim} it has been studied the following
impulsive process as partial sums in a series scheme
\begin{eqnarray}\label{1adf1}
\xi^{\varepsilon}(t)=\xi_0^{\varepsilon}+\sum_{k=1}^{\nu(t)}\alpha^{\varepsilon}_k(x^{\varepsilon}_{k-1}),\quad
t\ge 0,
\end{eqnarray}
the random variables $\alpha_k^{\varepsilon} (x), k \geq 1$ are
supposed to be independent and perturbed by the jump Markov process
$x(t), t\ge 0$.

We propose to study generalization of the problem (\ref{1adf1}):
\begin{eqnarray}\label{1bdf1}
\xi^{\varepsilon}(t)=\xi_0^{\varepsilon}+\sum_{k=1}^{\nu(t)}\alpha^{\varepsilon}_k(\xi^{\varepsilon}_{k-1},x^{\varepsilon}_{k-1}),\quad
t\ge 0.
\end{eqnarray}
Here the random variables $\alpha_k^{\varepsilon} (u, x), k \geq 1$
depend on the process $\xi^{\varepsilon}(t)$.

We propose to study convergence of (\ref{1bdf1}) using a combination
of two methods. The method proposed by R. Liptser in \cite{lip3},
based on semimartingales theory, is combined with a solution of
singular perturbation problem instead of ergodic theorem. So, the
method includes two steps.

In the first step we prove the relative compactness of the
semimartingales representation of the family $\xi^\eps$, $\eps>0$,
by proving the following two facts as proposed in Liptser
\cite{lip1}:
$$\lim\limits_{c\to \infty}\sup\limits_{\varepsilon\leq\varepsilon_0}
\mathbf{P}\{\sup\limits_{t\leq T}|\xi^{\varepsilon}(t)|>c\}=0,$$
known as the compact containment condition, and
$$\mathbf{E}|\xi^{\varepsilon}(t)-\xi^{\varepsilon}(s)|^2\le k
|t-s|,$$ for some positive constant $k$.

In the second step we prove convergence of two components of Markov
process $\xi^{\varepsilon}(t),
{\ae}_t^{\varepsilon}:={\ae}(t/\varepsilon^2)$
by using singular perturbation technique as presented in
\cite{korlim}.

Finally, we apply Theorem 6.3 from \cite{korlim}.

The paper is organized as follows. In Section 2 we present the
time-scaled impulsive process (\ref{1bdf1}) and the switching Markov
process. In the same section we present the main results of L\'{e}vy
approximation. In Section 3 we present the proof of the theorem.

\vsp
\section{Main results}
Let us consider the space $\matR^d$ endowed with a norm
$\abso{\cdot}$ ($d\ge 1$), and $(E,\mathcal{E})$, a {\it standard
phase space}, (i.e., $E$ is a Polish space and $\mathcal{E}$ its
Borel $\sigma$-algebra). For a vector $v\in \matR^d$ and a matrix
$c\in \matR^{d\times d}$ , $v^*$ and $c^*$ denote their transpose
respectively. Let $C_3(\matR^d)$ be a measure-determining class of
real-valued bounded functions, such that $g(u)/\abso{u}^2 \to 0$, as
$\abso{u}\to 0$ for $g\in C_3(\matR^d)$ (see \cite{jacod1,korlim}).

We introduce a family of random sequences $\alpha^{\varepsilon}_{k}
(x), k = 1, 2, ..., x \in E$, where $E$ is a non-empty set, indexed
by the small parameter $\varepsilon>0$. For any $\varepsilon > 0$,
and any sequence $z_k, k \geq 0$, of elements of $\mathbb{R}^d\times
E$, the random variables $\alpha_k^{\varepsilon} (z_{k-1}), k \geq
1$ are supposed to be independent. Let us denote by
$G_{u,x}^{\varepsilon}$ the distribution function of
$\alpha_k^{\varepsilon} (x)$, that is,
$$G_{u,x}^{\varepsilon}(dv) := P(\alpha_k^{\varepsilon} (u,x) \in
dv), k \geq 0, \varepsilon
> 0, x \in E, u\in \mathbb{R}^d.$$

The switching Markov process ${x}(t), t\ge 0$ on the standard phase
space $(E,\mathcal{E})$, is defined by the generator
\begin{eqnarray}\label{1gen0}
\mathbf{Q}\varphi(x) = q(x)\int_E P(x,dy)[\varphi(y)-\varphi(x)],
\end{eqnarray}
where $q(x), x\in E$, is the intensity of jumps function of $x(t),
t\ge 0$, and $P(x,dy)$ the transition kernel of the embedded Markov
chain $x_n,n \ge 0$, defined by $x_n= x(\tau_n), n\ge 0$, with
$0=\tau_0\le \tau_1\le ...\le \tau_n\le ...$ the jump times of
$x(t), t\ge 0$. Corresponding counting processes of jumps $\nu(t):=
\max\{k \geq 0 : \tau_k \leq t\}$. We make natural assumptions for
the counting process $\nu(t)$, namely:
\begin{eqnarray}\label{1con1}
\int_0^t\mathbf{E}[\varphi(s)d\nu(s)]<l_1\int_0^t\mathbf{E}(\varphi(s))ds
\end{eqnarray}
for any nonnegative, increasing $\varphi(s)$ and $l_1>0.$

Now we define a family of jump Markov processes
$x^{\varepsilon}(t):=x(t/\varepsilon^2), t \geq 0$, with embedded
Markov renewal process $x^{\varepsilon}_k, \tau^{\varepsilon}_k , k
\geq 0$, and counting processes of jumps
$\nu^{\varepsilon}(t)=\nu(t/\varepsilon^2), t \geq 0$. Thus, times
$\tau^{\varepsilon}_k , k \geq 0$, are jump times,
$x_k^{\varepsilon} := x^{\varepsilon}(\tau^{\varepsilon}_k )$, and
$\nu^{\varepsilon}(t) := \max\{k \geq 0 : \tau^{\varepsilon}_k \leq
t\}$.

The impulsive processes $\xi^{\varepsilon}(t), t\geq 0,
\varepsilon>0$ on $\mathbb{R}^d$ in the series scheme with small
series parameter $\varepsilon\to 0$, $(\varepsilon>0)$ are defined
by the sum (\cite[Section 9.2.1]{korlim})
\begin{eqnarray}\label{1adf2}
    \xi^{\varepsilon}(t)=\xi_0^{\varepsilon}+\sum_{k=1}^{\nu(t/\varepsilon^2)}\alpha^{\varepsilon}_{k}(\xi^{\varepsilon}_{k-1},x^{\varepsilon}_{k-1}),\quad
t\ge 0.
\end{eqnarray}

Here
$$\xi^{\varepsilon}_n:=\xi(\varepsilon^2\tau_n)=\xi_0^{\varepsilon}+\sum_{k=1}^n\alpha_k^{\varepsilon}(\xi^{\varepsilon}_{k-1},x^{\varepsilon}_{k-1}).$$

 It is worth noticing that the coupled process $\xi^\varepsilon(t),
x^{\varepsilon}(t), t \geq 0$, is a Markov additive process (see,
e.g., \cite[Section 2.5]{korlim}).

The L\'{e}vy approximation of Markov impulsive process (\ref{1adf2})
is considered under the following conditions.

\begin{description}
\item[C1:] The Markov process ${x}(t), t \geq 0$ is uniformly
ergodic with $\pi(B), B\in \mathcal{E}$ as stationary distribution.

\item[C2:] {\it L\'{e}vy approximation}. The family of impulsive processes $\xi^{\varepsilon}(t), t\geq 0$ satisfies the
L\'{e}vy approximation conditions \cite[Section 7.2.3]{korlim}.
\begin{description}
\item[L1:] Initial value condition
$$\sup\limits_{\varepsilon>0} E|\xi_0^{\varepsilon}|\leq C < \infty.$$

\item[L2:]Approximation of the mean values:
$$a^{\varepsilon}(u;x) = \int_{\mathbb{R}^d} vG^{\varepsilon}_{u,x}(dv)
= \varepsilon a_1(u;x)+\varepsilon^2 [a(u;x) +\theta_a^{\varepsilon}
(u;x)],$$ and
$$c^{\varepsilon}(u;x) = \int_{\mathbb{R}^d}
vv^*G^{\varepsilon}_{u,x}(dv) = \varepsilon^2 [c(u;x) +
\theta_c^{\varepsilon} (u;x)],$$ where functions $a_1, a$ and $c$
are bounded.

\item[L3:] Poisson approximation condition for intensity kernel (see \cite{jacod1})
$$G_g^{\varepsilon}(u;x) = \int_{\mathbb{R}^d} g(v)G^{\varepsilon}_{u,x}(dv)
= \varepsilon^2[G_g(u;x) + \theta^{\varepsilon}_g(u;x)]$$ for all $g
\in C_3(\mathbb{R}^d)$, and the kernel $G_g(u;x)$ is bounded for all
$g \in C_3(\mathbb{R}^d)$, that is,
$$|G_g(u;x)| \leq G_g \quad \hbox{(a constant depending on $g$)}.$$

Here \begin{eqnarray}\label{1gg2}
   G_{g}(u;x) =\int_{\mathbb{R}^d} g(v)G_{u,x}(dv),\quad g \in C_3(\mathbb{R}^d).
\end{eqnarray}

The above negligible terms $\theta_a^\eps,\theta_c^\eps,
\theta_g^\eps$ satisfy the condition $$\sup\limits_{x\in E}
|\theta_{\cdot}^{\varepsilon}(u;x)|\to 0,\quad  \varepsilon\to 0.$$

\item[L4:] {\it Balance condition}.
$$\int_E\rho(dx)a_1(u;x)=0.$$
\end{description}

In addition the following conditions are used:
\item[C3:] {\it Uniform square-integrability}:
$$\lim\limits_{c\to\infty}\sup\limits_{x\in E} \int_{|v|>c} vv^*G_{u,x}(dv) = 0,$$
where the kernel $G_{u,x}(dv)$ is defined on the measure determining
class $C_3(\mathbb{R}^d)$ by the relation (\ref{1gg2}).

\item[C4:] {\it Linear growth}: there exists a positive constant $L$ such that
$$|a(u;x)|\leq L(1+|u|),\quad\hbox{and}\quad |c(u;x)|\leq L(1+\abso{u}^2),$$
and for any real-valued non-negative function $f(x), x\in
\mathbb{R}^d$, such that $\int_{\mathbb{R}^d\setminus
\{0\}}(1+f(x))\abso{x}^2dx<\infty,$ we have
$$|G_{u,x} (v)|\leq Lf(v)(1+\abso{u}).$$
\end{description}

\vsp
The main result of our work is the following.

\begin{thm} Under conditions $\mathbf{C1-C4}$
the weak convergence
$$\xi^{\varepsilon}(t)\Rightarrow \xi^0(t),\quad \varepsilon \to 0$$ takes
place.

The limit process $\xi^0(t), t\geq0$ is a L\'{e}vy process defined
by the generator $\mathbf{L}$ as follows
\begin{eqnarray}\label{1limgen}
\mathbf{L}\varphi(u)=(\widehat{a}(u)-\widehat{a}_0(u))\varphi'(u)+
\frac{1}{2}\sigma^2(u)\varphi''(u) + \lambda(u)\int_{\mathbb{R}^d}
[\varphi(u + v)-\varphi(u)]G_{u}^0(dv),
\end{eqnarray}
with $\sigma^2(u)\geq 0$, where:
$$\widehat{a}(u)=q\int_E\rho(dx)a(u;x), \widehat{a}_0(u)=\int_EvG_u(dv), G_u(dv)=q\int_E\rho(dx)G_{u,x}(dv),$$
$$\sigma^2(u)=2\int_E\pi(dx)[\widetilde{a}_1(u;x)R_0\widetilde{a}_1^*(u;x)], \widetilde{a}_1(u;x):=q(x)\int_EP(x,dy)a_1(u;x)$$
\hskip40mm $\lambda(u)=qG_u(\mathbb{R}^d),$ \hskip5mm
$G_{u}^0(dv)=G_u(dv)/G_u(\mathbb{R}^d).$
\end{thm}


\vsp
\section{Proof of Theorem 1}
The proof of Theorem 1 is based on the semimartingale representation
of the impulsive process (\ref{1adf2}).

We split the proof of Theorem 1 in the following two steps.

\noindent {\sc Step 1}. In this step we establish the relative
compactness of the family of processes $\xi^{\varepsilon}(t), t\geq
0, \varepsilon>0$ by using the approach developed in \cite{lip3}.
Let us remind that the space of all probability measures defined on
the standard space $(E,\mathcal{E})$ is also a Polish space; so the
relative compactness and tightness are equivalent.

First we need the following lemma.

\begin{lm} Under assumption $\mathbf{C4}$ there exists a
constant $k>0$, independent of $\varepsilon$ and dependent on $T$,
such that
$$\mathbf{E}\sup\limits_{t\leq T}|\xi^{\varepsilon}(t)|^2\leq k_T.$$
\end{lm}

\begin{cor}  Under assumption $\mathbf{C4}$, the following
compact containment condition (CCC) holds:
$$\lim\limits_{c\to \infty}\sup\limits_{\varepsilon\leq\varepsilon_0}
\mathbf{P}\{\sup\limits_{t\leq T}|\xi^{\varepsilon}(t)|>c\}=0.$$
\end{cor}
\noindent{\it Proof}: The proof of this corollary follows from
Kolmogorov's inequality.\fin

\vsp \noindent{\it Proof of Lemma 1}: (following \cite{lip3}). The
impulsive process (\ref{1adf2}) has the following semimartingale
representation
\begin{eqnarray}\label{1smdecomp}
\xi^{\varepsilon}(t)=u+B_t^{\varepsilon}+M_t^{\varepsilon},
\end{eqnarray}
where $u= \xi^\eps_0$; $B_t^{\varepsilon}$ is the predictable drift
$$B_t^{\varepsilon}=\sum_{k=1}^{\nu(t/\varepsilon^2)}a^\eps(\xi^{\varepsilon}_{k-1},{x}^{\varepsilon}_{k-1})
=A_1^\varepsilon(t)+A^\varepsilon(t)+\theta^\varepsilon_a(t),$$
where
$$A^{\varepsilon}_1(t):=\varepsilon\sum_{k=1}^{\nu(t/\varepsilon^2)}a_{1}
(\xi_{k-1}^{\varepsilon}, x_{k-1}^{\varepsilon}),
A^{\varepsilon}(t):= \varepsilon^2\sum_{k=1}^{\nu(t/\varepsilon^2)}a
(\xi_{k-1}^{\varepsilon} ,x_{k-1}^{\varepsilon}).$$

\begin{eqnarray}\label{1qch}
\langle M^{\varepsilon}\rangle_t&=&\varepsilon^2\sum_{k=1}^{\nu
(t/\varepsilon^2)}
c(\xi^{\varepsilon}_{k-1};{x}_{k-1}^{\varepsilon})
+\varepsilon^2\sum_{k=1}^{\nu
(t/\varepsilon^2)}\int_{\mathbb{R}^d\setminus\{0\}}vv^*G(\xi^{\varepsilon}_{k-1},dv;{x}_{k-1}^{\varepsilon})+
\theta^{\varepsilon}_c(t),
\end{eqnarray} and for every finite $T>0$
$$\sup\limits_{0\leq t\leq T} |\theta^\varepsilon_{\cdot}(t)|\rightarrow 0, \varepsilon\rightarrow 0.$$

To verify compactness of the process $\xi^{\varepsilon}(t)$ we split
it at two parts.

The first part of order $\varepsilon$
$$A_1^\varepsilon(t)=\varepsilon\sum_{k=1}^{\nu(t/\varepsilon^2)}
a_1(\xi^{\varepsilon}_{k-1};{x}^{\varepsilon}_{k-1}),$$ can be
characterized by the martingale
$$\widetilde{\mu}_t^{\varepsilon}=\varphi^{\varepsilon}(A_1^\varepsilon(t))+\varphi^{\varepsilon}(A_1^\varepsilon(0))-
\int_0^t\mathbf{L}^{\varepsilon}\varphi^{\varepsilon}(A_1^\varepsilon(s))ds.$$

Thus (see, for example Theorem 1.2 in \cite{korlim}), it has
quadratic characteristic
$$<\widetilde{\mu}^{\varepsilon}>_t=\int_0^t\left[\mathbf{L}^{\varepsilon}(\varphi^{\varepsilon}
(A_1^\varepsilon(s)))^2-2\varphi^{\varepsilon}(A_1^\varepsilon(s))
\mathbf{L}^{\varepsilon}\varphi^{\varepsilon}(A_1^\varepsilon(s))\right]ds.$$

Applying the operator
$\mathbf{L}^{\varepsilon}=\varepsilon^{-2}\mathbf{Q}+\varepsilon^{-1}\mathbf{A}_1$
to test-function
$\varphi^{\varepsilon}=\varphi+\varepsilon\varphi_1$ (here
$\mathbf{A}_1(u;x)\varphi(v)=a_1(u;x)\varphi'(v)$) we obtain the
integrand of the view
$$Q\varphi_1^2-2\varphi_1Q\varphi_1.$$

It is independent of $\varepsilon$ and limited. The boundedness of
the quadratic characteristic provides
$\widetilde{\mu}_t^{\varepsilon}$ is compact. Thus,
$\varphi(A_1^\varepsilon(t))$ is compact too and bounded uniformly
by $\varepsilon$.

Now we should study the second part of order $\varepsilon^2$.

For a process $y(t), t\ge 0$, let us define the process
$y^\dag(t)=\sup\limits_{s\leq t}|y(s)|,$ then from (\ref{1smdecomp})
we have
\begin{eqnarray}\label{1eq4}
((\xi^{\varepsilon}(t))^\dag)^2\le
4[u^2+((A^{\varepsilon}(t))^\dag)^2+((M^{\varepsilon}_t)^\dag)^2].
\end{eqnarray}

Now we may apply the result of Section 2.3 \cite{korlim}, namely
$$\sum_{k=1}^{\nu(t)} a(\xi^\eps_{k-1}, x^\eps_{k-1})=\int_0^ta(\xi^{\varepsilon}(s),x^{\varepsilon}(s))d\nu(s).$$

Condition $\mathbf{C4}$  implies that for sufficiently large
$\varepsilon$
\begin{eqnarray}\label{1eq51}
(A^\varepsilon(t))^\dag && =
\varepsilon^2\int_0^{t/\varepsilon^2}a(\xi^{\varepsilon}(s),x^{\varepsilon}(s))d\nu(s)\leq
L\varepsilon^2\int_0^{t/\varepsilon^2}(1+(\xi^{\varepsilon}(s))^\dag)d\nu(s)\end{eqnarray}

Now, by Doob's inequality (see, e.g., \cite[Theorem 1.9.2]{lip1}),
$$\mathbf{E}((M_t^{\varepsilon})^\dag)^2\leq
4\abso{\mathbf{E}\langle M^{\varepsilon}\rangle_t},$$ (\ref{1qch})
and condition \textbf{C4} we obtain
\begin{eqnarray}\label{1eq6}
\abso{\langle
M^{\varepsilon}\rangle_t}=\abso{\varepsilon^2\int_0^{t/\varepsilon^2}
c(\xi^{\varepsilon}(s);{x}^{\varepsilon}(s))
+\varepsilon^2\int_0^{t/\varepsilon^2}\int_{\mathbb{R}^d\setminus
\{0\}}vv^*G(\xi^{\varepsilon}(s),dv;{x}_s^{\varepsilon})d\nu(s)}\leq\nonumber\\
2L(1+r_1)\varepsilon^2\int_0^{t/\varepsilon^2}[1+((\xi^{\varepsilon}(s))^\dag)^2]d\nu(s),
\end{eqnarray}where
$r_1=\int_{\mathbb{R}^d\setminus \{0\}}\abso{x}^2f(x)dx.$

Inequalities (\ref{1eq4})-(\ref{1eq6}), condition (\ref{1con1}) and
Cauchy-Bunyakovsky-Schwarz inequality,
([$\int_0^t\varphi(s)ds]^2\leq t\int_0^t\varphi^2(s)ds$), imply
$$\mathbf{E}((\xi^{\varepsilon}(t))^\dag)^2\leq
k_1+k_2\varepsilon^2\int_0^{t/\varepsilon^2}\mathbf{E}[((\xi^{\varepsilon}(s))^\dag)^2d\nu(s)]\leq
k_1+k_2l_1\varepsilon^2\int_0^{t/\varepsilon^2}\mathbf{E}((\xi^{\varepsilon}(s))^\dag)^2ds=$$
$$ k_1+k_2l_1\int_0^{t}\mathbf{E}((\xi^{\varepsilon}(s))^\dag)^2ds,$$
where $k_1, k_2$ and $l_1$ are positive constants independent of
$\varepsilon$.

By Gronwall inequality (see, e.g., \cite[p. 498]{ethier}), we obtain
$$\mathbf{E}((\xi^{\varepsilon}(t))^\dag)^2\leq k_1\exp(k_2l_1 t).$$

Thus, both parts of $\xi^{\varepsilon}(t)$ are compact and bounded,
so $$\mathbf{E}\sup\limits_{t\leq T}|\xi^{\varepsilon}(t)|^2\leq
k_T.$$

Hence the lemma is proved. \fin

\begin{lm} Under assumption $\mathbf{C4}$ there exists a
constant $k>0$, independent of $\varepsilon$ such that
$$\mathbf{E}|\xi^{\varepsilon}(t)-\xi^{\varepsilon}(s)|^2\leq k |t-s|.$$
\end{lm}

\noindent{\it Proof}: In the same manner with (\ref{1eq4}), we may
write
$$|\xi^{\varepsilon}(t)-\xi^{\varepsilon}(s)|^2\leq 2|B_t^{\varepsilon}
-B_s^{\varepsilon}|^2+2|M_t^{\varepsilon}-M_s^{\varepsilon}|^2.$$ By
using Doob's inequality, we obtain
$$\mathbf{E}|\xi^{\varepsilon}(t)-\xi^{\varepsilon}(s)|^2\leq
2\mathbf{E}\{|B_t^{\varepsilon}-B_s^{\varepsilon}|^2+8\abso{\langle
M^{\varepsilon}\rangle_t-\langle M^{\varepsilon}\rangle_s}\}.$$

Now (\ref{1eq6}) and condition (\ref{1con1}) and assumption
$\mathbf{C4}$ imply
$$|B_t^{\varepsilon}-B_s^{\varepsilon}|^2+8\abso{\langle
M^{\varepsilon}\rangle_t-\langle M^{\varepsilon}\rangle_s}\leq
k_3[1+((\xi^{\varepsilon}(T))^\dag)^2]|t-s|,$$ where $k_3$ is a
positive constant independent of $\varepsilon$.

From the last inequality and Lemma 1 the desired conclusion is
obtained. \fin

The conditions proved in Corollary 2 and Lemma 2 are necessary and
sufficient for the compactness of the family of processes
$\xi^{\varepsilon}(t), t\geq 0, \varepsilon>0$.

\vsp \noindent{\sc Step 2}. At the next step of proof we apply the
problem of singular perturbation to the generator of the process
$\xi^{\varepsilon}(t).$ To do this, we mention the following
theorem.
$C^2_0(\mathbb{R}^d\times E)$ is the space of real-valued twice
continuously differentiable functions on the first argument, defined
on $\mathbb{R}^d\times E$ and vanishing at infinity, and
$C(\mathbb{R}^d\times E)$ is the space of real-valued continuous
bounded functions defined on $\mathbb{R}^d\times E$.

\begin{thm}(\cite[Theorem 6.3]{korlim}) Let the following conditions hold
for a family of Markov processes $\xi^{\varepsilon}(t), t\ge 0,
\varepsilon>0$:
\begin{description}
\item[CD1:] There exists a family of test functions
$\varphi^{\varepsilon}(u, x)$ in $C^2_0(\mathbb{R}^d\times E)$, such
that $$\lim\limits_{\varepsilon\to 0}\varphi^{\varepsilon}(u, x) =
\varphi(u),$$ uniformly on $u, x.$

\item[CD2:] The following convergence holds
$$\lim\limits_{\varepsilon\to
0}\mathbf{L}^{\varepsilon}\varphi^{\varepsilon}(u, x) =
\mathbf{L}\varphi(u),$$ uniformly on $u, x$. The family of functions
$\mathbf{L}^{\varepsilon}\varphi^{\varepsilon}, \varepsilon>0$ is
uniformly bounded, and $\mathbf{L}\varphi(u)$ and
$\mathbf{L}^{\varepsilon}\varphi^{\varepsilon}$ belong to
$C(\mathbb{R}^d\times E)$.

\item[CD3:] The quadratic characteristics of the
martingales that characterize a coupled Markov process
$\xi^{\varepsilon}(t), x^{\varepsilon}(t), t\geq0, \varepsilon>0$
have the representation $\left\langle
\mu^{\varepsilon}\right\rangle_t = \int^t_0
\zeta^{\varepsilon}(s)ds,$ where the random functions
$\zeta^{\varepsilon}, \varepsilon> 0,$ satisfy the condition
$$\sup\limits_{0\leq s \leq T} \mathbf{E}|\zeta^{\varepsilon}(s)|\leq
c < +\infty.$$

\item[CD4:] The convergence of the initial values holds
and
$$\sup\limits_{\varepsilon>0}\mathbf{E}|\zeta^{\varepsilon}(0)|\leq C
< +\infty.$$
\end{description}

Then the weak convergence
$$\xi^{\varepsilon}(t)\Rightarrow \xi(t),\quad \varepsilon\to 0,$$
takes place.
\end{thm}

We consider the two component Markov process
$\xi^{\varepsilon}(t),{x}_t^{\eps}, t\ge 0$ which can be
characterized by the martingale
$$\mu_t^{\varepsilon}=\varphi(\xi^{\varepsilon}(t),{x}_t^{\eps})-
\int_0^t\mathbf{L}^{\varepsilon}
\varphi(\xi^{\varepsilon}(s),{x}_t^{\eps})ds,$$ where its generator
$\mathbf{L}^{\varepsilon}$ has the following representation
\cite{korlim} Lemma 9.1
\begin{eqnarray}\label{1eq7}
\mathbf{L}^{\varepsilon}\varphi(u,x)=\varepsilon^{-2}q(x)\left[
\int_E P(x, dy) \int_{\mathbb{R}^d}
G_{u,y}^{\varepsilon}(dz)\varphi(u + z,y)-\right.\\\nonumber
\left.\varphi(u,x)\right].
\end{eqnarray}

By analogy with \cite[Lemma 9.2]{korlim} we may prove the following
result:

\begin{lm} The main part in the asymptotic representation of the generator (\ref{1eq7}) is as follows
$$\mathbf{L}^{\varepsilon}\varphi(u,v,x) = \varepsilon^{-2}\mathbf{Q}\varphi(\cdot,\cdot,x) +
\varepsilon^{-1}\mathbf{Q}_0a_1(u;x)\varphi'_u(u,\cdot,\cdot) +
\mathbf{Q}_0[a(u;x) - a_0(u;x)]\varphi'_u(u,\cdot,\cdot) +$$ $$
\mathbf{Q}_0\mathbf{G}_{u,x}\varphi(u,\cdot,\cdot)
$$
where:
$$\mathbf{Q}_0\varphi(x) := q(x) \int_E P(x, dy)\varphi(y),
\mathbf{G}_{u,x}\varphi(u) := \int_{\mathbb{R}^d} [\varphi(u + v)
-\varphi(u)]G_{u,x}(dv),$$ $$a_0(u;x)=\int_EvG_{u,x}(dv).$$
\end{lm}

\noindent{\it Proof} of this Lemma is analogical to the proof of
\cite[Lemma 9.2]{korlim}.

The solution of the singular perturbation problem at the test
functions
$\varphi^{\varepsilon}(u,x)=\varphi(u)+\varepsilon\varphi_1(u,x)+\varepsilon^2\varphi_2(u,x)$
in the form $\mathbf{L}^{\varepsilon}\varphi^{\varepsilon}
={\mathbf{L}}\varphi+\theta^{\varepsilon}\varphi$ can be found in
the same manner with Lemma 9.3 in \cite{korlim}. That is
\begin{eqnarray}\label{1eq8}
{\mathbf{L}}=\Pi[\mathbf{Q}_0(\mathbf{A}(x)+\mathbf{G}_{u,x})+\mathbf{Q}_0\mathbf{A}_1(x)R_0\mathbf{Q}_0\mathbf{A}_1(x)]\Pi,
\end{eqnarray} where $$\mathbf{A}(x)\varphi(u):=[a(u;x)-a_0(u;x)]\varphi'(u), \mathbf{A}_1(x)\varphi(u):=a_1(u;x)\varphi'(u).$$

Simple calculations give us (\ref{1limgen}) from (\ref{1eq8}).

Now Theorem 2 can be applied.

We see from (\ref{1eq7}) and (\ref{1eq8}) that the solution of
singular perturbation problem for
$\mathbf{L}^{\varepsilon}\varphi^{\varepsilon}(u,v;x)$ satisfies the
conditions \textbf{CD1, CD2}. Condition \textbf{CD3} of this theorem
implies that the quadratic characteristics of the martingale,
corresponding to a coupled Markov process, is relatively compact.
The same result follows from the CCC (see Corollary 2 and Lemma 2)
by \cite{jacod1}. Thus, the condition \textbf{CD3} follows from the
Corollary 2 and Lemma 2. As soon as $\xi^{\varepsilon}(0)=\xi^0(0)$
we see that the condition \textbf{CD4} is also satisfied. Thus, all
the conditions of above Theorem 2 are satisfied, so the weak
convergence $\xi^{\varepsilon}(t)\rightarrow \xi^0(t)$ takes place.

\vsp

Theorem 1 is proved.\fin


\begin{thebibliography}{1}

\bibitem{ber} Bertoin J. (1996). {\it L\'{e}vy processes.} Cambridge Tracts in Mathematics,
121. Cambridge University Press, Cambridge.


\bibitem{ethier} Ethier S.N., Kurtz T.G. (1986). {\it Markov Processes:
Characterization and convergence}, J. Wiley, New York.

\bibitem{gisk} Gihman, I.I., Skorohod, A.V. (1974). {\it Theory of stochastic processes,
vol. 1,2,3,} Springer, Berlin.

\bibitem{jacod1} Jacod J., Shiryaev A.N. (1987). {\it Limit
Theorems for Stochastic Processes}, Springer-Verlang, Berlin.

\bibitem{korlim} Koroliuk V.S., Limnios N. (2005). {\it Stochastic Systems in Merging Phase
Space}, World Scientific, Singapore.


\bibitem{lip3} Liptser R. Sh. (1994). The Bogolubov averaging principle for
semimartingales, {\it Proceedings of the Steklov Institute of
Mathematics}, Moscow, No 4, 12 pages.

\bibitem{lip1} Liptser R. Sh., Shiryayev A. N. (1989). {\it Theory of
Martingales}, Kluwer Academic Publishers, Dordrecht, The
Netherlands.

\bibitem{sat} Sato K.-I. (1999). {\it L\'{e}vy processes and infinitely divisible
distributions.} Cambridge Studies in Advanced Mathematics, 68.
Cambridge University Press, Cambridge.


\end{thebibliography}
\end{document}